# Maximum Likelihood Estimations Based on Upper Record Values for Probability Density Function and Cumulative Distribution Function in Exponential Family and Investigating Some of Their Properties

## S. D. Gore[1], S. Hosseini[2], P. Nasiri[3]

**Abstract**

In this paper a useful subfamily of the exponential family has been considered. The ML estimation based on upper record values has been calculated for the parameter, Cumulative Density Function, and Probability Density Function of the family. Also, the relations between MLE based on record values and a random sample has been discussed. Additionally, some properties of these estimators has been investigated. Finally, it has been proven that these estimators have some useful properties for samples with large size.

**1. Introduction:**

An exponential family includes a wide range of statistical distributions in two discrete and continuous states, that has a large importance in the distribution theory, and with which it might be possible to integrally investigate many properties of the distribution. Different forms have been presented for this family in the one-parametric state. Generally, a one-variable exponential family is a group of distributions that their probability density function is as follows:

$$f(x;\theta) = h(x)g(\theta)\exp\{-\eta(\theta)T(x)\}$$

Here, it is considered an important subfamilies from the exponential family.

$$F(x;\theta) = 1 - \exp\{-B(\theta)A(x)\},$$

in which A is an increasing function and

$$B(\theta) > 0; x \in [a,b]; A(a) = 0, A(b) = +\infty; a,b \in R$$

It is called the first type exponential family for ease in this paper. Many one-parametric continuous distributions can be studied by choosing the above general form. The following tables describe this subject:

---
[1] Professor at department of statistics, Savirtibai Phule Pune University, India
[2] Correspondent author and lecturer at department of accounting, Cihan University of Erbil, Iraq. Email: S.hosseini.stat@gmail.com
[3] Associate professor at Payam Noor University of Tehran

**Table 1: Instances for the exponential family of the first type**

| $A(x)$ | $B(\theta)$ | Distribution Name |
|---|---|---|
| $A(x) = x$ | $\theta$ | Exponential |
| $A(x) = Log(1+x)$ | $B(\theta) = \dfrac{1}{\theta}$ | Lomax |
| $A(x) = x^{\alpha}$ | $B(\theta) = \theta$ | Weibull |
| $A(x) = -Ln(x)$ | $B(\theta) = \theta$ | Pareto |

Other distributions, which belong to the mentioned family, may be found by more investigations. Therefore, the reason for choosing the above general family lies in its wideness and comprehensiveness.

Probability density function (PDF) diagram for the distributions in Table.1 and their cumulative distribution function (CDF) diagram are shown in Figures 1.1 and 1.2, respectively.

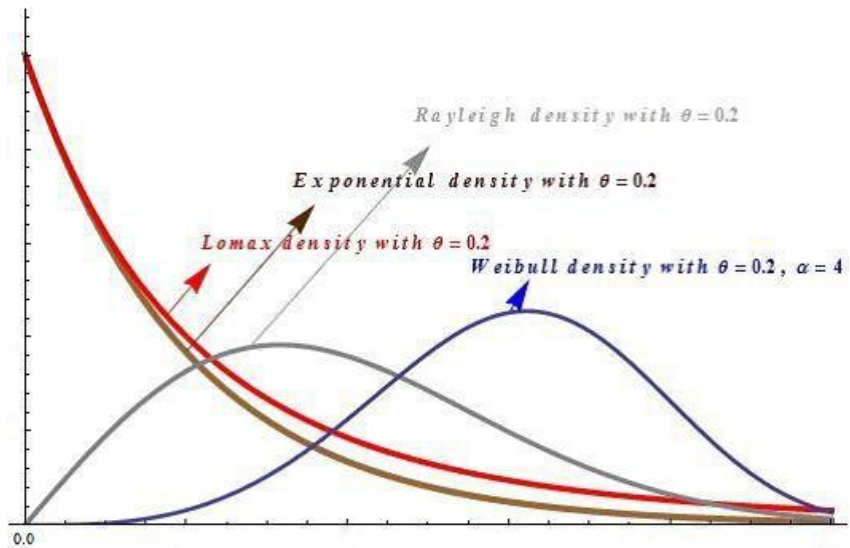

Figure 1.1

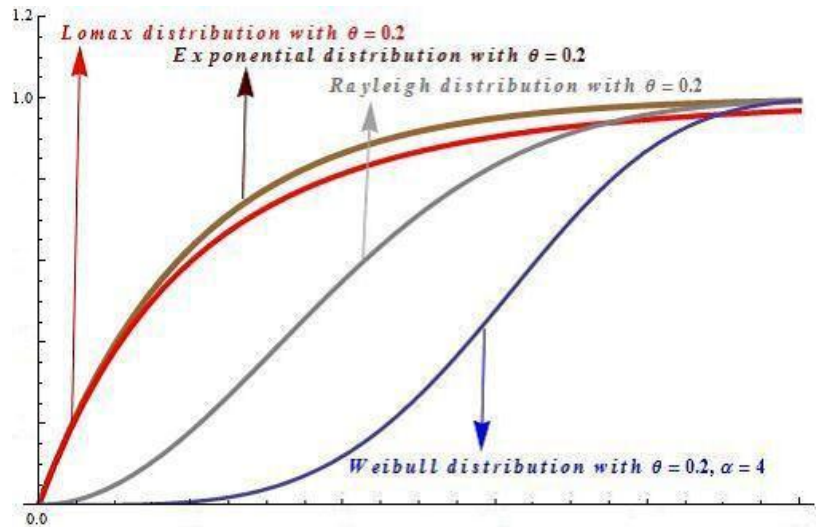

Figure 1.2

Determination and investigation of the properties of the maximum likelihood estimations based on record values for probability density function and cumulative distribution function in the given exponential family is the goal of this paper; thus the maximum likelihood estimations for the unknown parameters of the exponential family, probability density function, and cumulative distribution function of the exponential family have been generally obtained based on a random sample and also values of upper records. In addition, validity of maximum likelihood estimators based on record values in different states has been presented along with a series of examples. Likewise, the criterion of Mean Square Error (MSE) has been considered in order to investigate and compare them. Additionally, the properties of these estimators in asymptotic conditions are interested by this paper. Chandler suggested for the first time the concept of records under a mathematical model (Chandler, 1952). During a very short period after him, this concept was changed and completed by many statisticians and the theory of records was formed which has very wide applications in analysis of competitive phenomena such as commerce and economics. After 21 years from the introduction of Chandler to this discussion, a wide range of researches and other works were performed in this field, and some interesting conclusions were resulted. The main conclusion was offered by Shorrock in 1973 when completing the asymptotic theory of records (Shorrock 1973). In addition, Ahsanullah in 1973 covered the main concepts of the theory of records in 1998 in a paper named "an introduction to record values". The theory of records, nowadays, is utilized in various statistical issues like the estimation theory, queuing theory, and etc. Among the works recently done in the field of estimation of distribution parameters based on records, it might be pointed out to those of Ahmadi & Doostparast (2006), Ahmadi, Ja'fari Jozani, Eric Marchand & Parsian (2009), Baklizi (2008) and Ragheb (2007), Nasiri and Hosseini (2012). Ahmadi & Balakrishnan in a paper published in 2004, did investigate the interval estimations based on record values for quantiles. As well, Nasiri and Hosseini in 2013 used a statistic based on record values to estimate Bayesian estimation for exponential distribution.

## 2. Maximum Likelihood Estimations Based on Random Sample and Upper Record Values for Unknown Parameter of Exponential Family

Considering an iid random sample of size n from the exponential family of the first type, the joint probability density function will be obtained as following:

$$f(x_1,...,x_n;\theta) = B^n(\theta)\prod_{i=1}^{n} A'(x_i)\exp\{-B(\theta)\sum_{i=1}^{n} A(x_i)\} \qquad (1)$$

The following fact has been used in the above relation:

$$f(x;\theta) = \frac{dF(x;\theta)}{dx} = A'(x)B(\theta)\exp\{-A(x)B(\theta)\}$$

Hence, the maximum likelihood function is determined as follows considering the joint probability density function $(X_1,...,X_n)$ (1).

$$L = \log f(x_1,...,x_n;\theta) = n\log B(\theta) + \sum_{i=1}^{n}\log A'(X_i) - B(\theta)\sum_{i=1}^{n} A(X_i)$$

On the other hand, it is known that the maximum likelihood estimations are determined from solving the equation $\frac{\partial L}{\partial \theta} = 0$, therefore:

$$\frac{\partial L}{\partial \theta} = 0 \Rightarrow \frac{nB'(\theta)}{B(\theta)} - B'(\theta)\sum_{i=1}^{n} A(X_i) = 0 \Rightarrow B(\theta) = \frac{n}{\sum_{i=1}^{n} A(X_i)}.$$

Now, if $B(\theta)$ is a one-to-one function from $\Theta$ to R, the following expression is simply resulted:

$$\theta = B^{-1}(\frac{n}{\sum_{i=1}^{n} A(X_i)}) \qquad (2)$$

Since the probability density function $f(x;\theta) = A'(x)B(\theta)\exp\{-A(x)B(\theta)\}$ with the properties and conditions described for A and B is a continuous function, it is possible to obtain its first and second order derivatives.

$$\frac{\partial^2 L}{\partial \theta^2} = \frac{nB''(\theta)B(\theta) - n[B'(\theta)]^2}{B^2(\theta)} - B''(\theta)t, (t = \sum_{i=1}^{n} A(x_i)) \qquad (3)$$

In addition, it is simply observed by substituting (2) in (3) that:

$$\frac{\partial^2 L}{\partial \theta^2}\bigg|_{B(\theta)=\frac{n}{\sum_{i=1}^{n}A(X_i)}} = \frac{nB''(\theta)(\frac{n}{t}) - n[B'(\theta)]^2}{(\frac{n}{t})^2} - B''(\theta)t = -\frac{[B'(\theta)t]^2}{n} < 0$$

Considering the negative sign of the likelihood function, it is found out that $\theta = B^{-1}(\frac{n}{\sum_{i=1}^{n}A(X_i)})$ is a likelihood value. With respect to the described issues, ML estimation of the parameter $\theta$ based on the random sample $(X_1,...,X_n)$ will be attained as following:

$$\hat{\theta}_{MLE,(X_1,...,X_n)} = B^{-1}(\frac{n}{\sum_{i=1}^{n}A(X_i)}) \tag{4}$$

On the other hand, $A(X)$ is distributed as exponential distribution by the parameter $\frac{1}{B(\theta)}$ ( $A(X) \stackrel{Distribution}{=} Exp(\frac{1}{B(\theta)})$ ), because

$$P(A(X) \leq t) = P(X \leq A^{-1}(t)) = 1 - \exp\{-B(\theta)t\}$$

Besides, since $X_i{}_{i:1,...,n}$ is a random sample, $A(X_i){}_{i:1,...,n}$ is also a random sample and

$$A(X_i){}_{i:1,...,n} \stackrel{Distribution}{=} Exp(\frac{1}{B(\theta)}), \tag{5}$$

so

$$T = \sum_{i=1}^{n} A(X_i) \stackrel{Distribution}{=} Gamma(n, \frac{1}{B(\theta)}).$$

Therefore, regarding equations (4) and (5):

$$\hat{\theta}_{MLE,(X_1,...,X_n)} \stackrel{Distribution}{=} B^{-1}(\frac{n}{T}) \tag{6}$$

in which $T \stackrel{Distribution}{=} Gamma(n, \frac{1}{B(\theta)})$.

Considering a random sample of size m which has n values of the upper records $(n < m)$, the joint probability density function of the upper records for this sample is determined as follows (Arnold, Balakrishnan and Nagaraja, 1998):

$$f(r_1,\ldots,r_n;\theta) = f(r_n;\theta)\prod_{i=1}^{n-1} h(r_i;\theta) = B^n(\theta)\{\prod_{i=1}^{n} A'(r_i)\}\exp\{-B(\theta)+A(r_n)\}$$

As a result, the likelihood function based on upper records is obtained like follows:

$$L_{Records} = Log f(r_1,\ldots,r_n;\theta) = nLogB(\theta) + \sum_{i=1}^{n} LogA'(r_i) - B(\theta)A(r_n).$$

Maximum likelihood estimation based on upper record values for parameter $\theta$ will be obtained from solving the equation $\frac{\delta L_{Records}}{\delta \theta} = 0$; hence:

$$\frac{\delta L_{Records}}{\delta \theta} = 0 \Rightarrow \frac{nB'(\theta)}{B(\theta)} - B'(\theta)A(r_n) = 0 \Rightarrow \hat{\theta} = B^{-1}(\frac{n}{A(R_n)}). \tag{7}$$

On the other hand:

$$\frac{\partial^2 L_{Records}}{\partial \theta^2} = \frac{nB''(\theta)B(\theta) - n[B'(\theta)]^2}{B^2(\theta)} - B''(\theta)l, (l = A(r_n)) \tag{8}$$

Substituting (7) in (8) the following relation is obtained:

$$\frac{\partial^2 L_{Records}}{\partial \theta^2}\bigg|_{B(\theta)=\frac{n}{A(r_n)}} = \frac{nB''(\theta)(\frac{n}{l}) - n[B'(\theta)]^2}{(\frac{n}{l})^2} - B''(\theta)l = -\frac{[B'(\theta)l]^2}{n} < 0.$$

Which implies that (7) is a maximum likelihood estimation based on upper record values because it maximizes the likelihood function based on record values. As a result:

$$\hat{\theta}_{MLE,(R_1,\ldots,R_n)} = B^{-1}(\frac{n}{A(R_n)}). \tag{9}$$

On the other hand, since $A(R_i)_{i:1,\ldots,n}$ is distributed as gamma by the parameter $(i,\frac{1}{B(\theta)})$, then $A(R_n)$ is distributed as gamma by the parameter $(n,\frac{1}{B(\theta)})$ (Arnold, Balakrishnan and Nagaraja, 1998).

$$T = A(R_n) \overset{Distribution}{=} Gamma(n,\frac{1}{B(\theta)}). \tag{10}$$

Considering (9) and (10) simultaneously leads to:

$$\hat{\theta}_{MLE,(R_1,...,R_n)} \overset{Distribution}{=} B^{-1}(\frac{n}{T}) \qquad (11)$$

in which $T = A(R_n) \overset{Distribution}{=} Gamma(n, \frac{1}{B(\theta)})$.

The above-mentioned issues have been summarized in the following theorem:

**Theorem 1**

In the exponential family of the first type, the maximum likelihood estimator based on a random sample of size n and the maximum likelihood estimator based on the first n values of the upper records are identically distributed. In other words, if g is a real function then:

$$i - g(\hat{\theta}_{MLE,((X_1,...,X_n))}) \overset{Distribution}{=} g(\hat{\theta}_{MLE,(R_1,...,R_n)})$$

$$ii - MSE(g(\hat{\theta}_{MLE,((X_1,...,X_n))})) = MSE(g(\hat{\theta}_{MLE,(R_1,...,R_n)}))$$

**Proof:**

Regarding the described issues and the relations (6) and (11), the proof is obvious.

The second part of the theorem 1 can be regarded as the most primary reason for using record estimators. Namely, in likelihood estimation of any function of $\theta$, a random sample of size n and n values of the upper records have the same errors (if MSE is regarded as the error criterion). Some examples will be presented in the next sections in which the record estimations are more justifiable. Before this subject, however, it is useful to define the symbol $\alpha(n)$ as follows:

$$MSE(g(\hat{\theta}_{MLE,(X_1,...,X_n)})) = MSE(g(\hat{\theta}_{MLE,(R_1,...,R_n)})) =$$
$$E[\{g(B^{-1}(\frac{n}{T})) - g(\theta)\}^2] = E[g^2(B^{-1}(\frac{n}{T}))] - 2g(\theta)E[g(B^{-1}(\frac{n}{T}))] + (g(\theta))^2 = \alpha(n) \qquad (12)$$

**Result 1**

Considering the first part of theorem 1 and relation (12), if a number of m upper record values $R_1,...,R_m$ exist in a random sample $X_1,...,X_n$, then the MSE value of ML estimation based on record values in this sample equals with MSE of ML estimation based on a random sample of size m, i.e.

$$MSE(\hat{\theta}_{MLE,(R_1,...,R_m)}) = MSE(\hat{\theta}_{MLE-based-on-random-sample-of-size-m-(X_1,...,X_m)}) = \alpha(m) \qquad (13)$$

**Result 2**

The relation (12) and the result (1) provide us a new tool and idea in order to decide by using it that in which conditions employing the record estimations is better. Two following instances clarify this issue.

**Example 1.** If $A(x) = x$ and $B(\theta) = \dfrac{1}{\theta}$, then $B^{-1}(\dfrac{n}{T}) = \dfrac{T}{n}$; it is known that $B^{-1}(\dfrac{n}{T}) = \dfrac{T}{n}$ is distributed as gamma by the parameters $(n, \dfrac{\theta}{n})$ $(B^{-1}(\dfrac{n}{T}) = \dfrac{T}{n} \overset{Distribution}{=} Gamma(n, \dfrac{\theta}{n}))$. Now, considering the relation (12)

$$E[(B^{-1}(\dfrac{n}{T}))^2] = \theta^2 + \dfrac{\theta^2}{n}, \quad E[B^{-1}(\dfrac{n}{T})] = \theta$$

therefore

$$MSE = \alpha(n) = \dfrac{\theta^2}{n}.$$

With respect to the figure, it is readily concluded that for a random sample of size n, if m records of the upper records are available $(m < n)$, then the maximum likelihood estimation based on record values have always more errors.

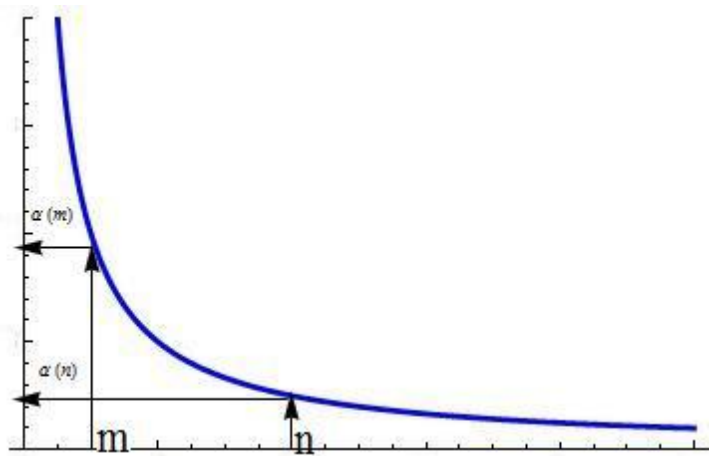

Figure 1.2. Comparison of maximum likelihood estimation
based on random sample and records in an exponential distribution.

As a consequence of Example 1, the maximum likelihood estimations based on record values might be inefficient and unsuitable in the first looking. It is shown that this belief is not in general the case and there are circumstances in which estimations based on record values are more efficient (i.e. have less MSEs). The following example clarifies this subject.

**Example 2.** Assuming that $A(x) = Log(\frac{x}{k})$ and $B(\theta) = \theta$, it is wanted to estimate the parameter $g(\theta) = k^\theta$.

As a consequence of the relation (12),

$$B^{-1}(\theta) = \theta, B^{-1}(\frac{n}{T}) = \frac{n}{T}, g(B^{-1}(\frac{n}{T})) = k^{\frac{n}{T}}$$

Thus

$$k^{\frac{n}{T}} = \exp\{Ln(k^{\frac{n}{T}})\} = \sum_{i=0}^{+\infty} \frac{(Ln(k^{\frac{n}{T}}))^i}{i!}, \text{ if for example } k = e \text{ then}$$

$$k^{\frac{n}{T}} = \exp\{Ln(k^{\frac{n}{T}})\} = \sum_{i=0}^{+\infty} \frac{(Ln(k^{\frac{n}{T}}))^i}{i!} = \sum_{i=0}^{+\infty} \frac{(\frac{n}{T})^i}{i!}.$$

Considering $E[\{g(B^{-1}(\frac{n}{T}))\}^\alpha] = E[\exp\{\frac{\alpha n}{T}\}] = E[\sum_{i=0}^{+\infty} \frac{(\frac{\alpha n}{T})^i}{i!}], \text{ and } \frac{(\frac{\alpha n}{T})^i}{i!} > 0$

Each variable $\frac{(\frac{\alpha n}{T})^i}{i!}$ is greater than zero then, that's why it is possible to interchange the integration and summation.

$$E[\sum_{i=0}^{+\infty} \frac{(\frac{\alpha n}{T})^i}{i!}] = \sum_{i=0}^{+\infty} E[\frac{(\frac{\alpha n}{T})^i}{i!}] = \sum_{i=0}^{+\infty} \{\int_0^{+\infty} \frac{\alpha^i n^i t^{n-1} \theta^n \exp\{-\theta t\}}{t^i i!(n-1)!} dt\} = \sum_{i=0}^{+\infty} \frac{\alpha^i n^i \theta^n}{i!(n-1)!} \{\int_0^{+\infty} t^{n-i-1} \exp(-\theta t) dt,$$

n and i are round numbers that's why $\int_0^{+\infty} t^{n-i-1} \exp(-\theta t) dt = \frac{\Gamma(n-i)}{\theta^{n-i}}.$  (14)

$$E[\{g(B^{-1}(\frac{n}{T}))\}^\alpha] = E[\sum_{i=0}^{+\infty} \frac{(\frac{\alpha n}{T})^i}{i!}] = \sum_{i=0}^{+\infty} \frac{\alpha^i n^i \theta^n}{i!(n-1)!} \{\int_0^{+\infty} t^{n-i-1} \exp(-\theta t) dt = \sum_{i=0}^{+\infty} \frac{\alpha^i n^i \theta^n}{i!(n-1)!} \frac{\Gamma(n-i)}{\theta^{n-i}}.$$

Having considered the above relation and the fact that Gamma function is defined for positive values, the MSE is readily resulted as follows:

$$MSE(g(\theta)) = \sum_{i=0}^{n-1} \frac{2^i n^i \theta^n}{i!(n-1)!} \frac{\Gamma(n-i)}{\theta^{n-i}} - 2e^\theta \sum_{i=0}^{n-1} \frac{n^i \theta^n}{i!(n-1)!} \frac{\Gamma(n-i)}{\theta^{n-i}} + e^{2\theta}.$$

If we draw $MSE(g(\theta))$ for different values of n, then the following figure is obtained:

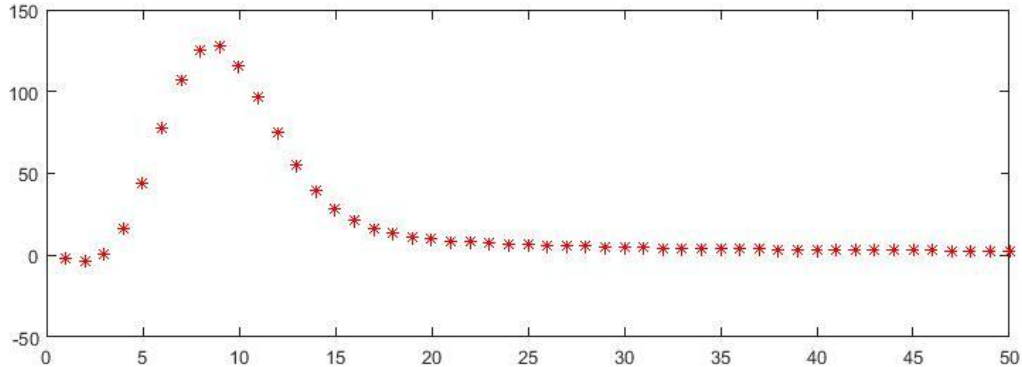

Figure 2.2. Different values of MSE in terms of n when estimating g

Consider the figure 2.2; if for instance the number of upper records in a sample of size 10 or 11 equals to 5 or 6, then the maximum likelihood estimation $g(\theta)$ based on records does have less amount of error. For another example, it is seen, with respect to the figure 2.2, that $\alpha(4) < \alpha(5) < \alpha(7) < \alpha(12)$, that is, if a sample of size 12 has 4 values of the upper records, then the ML estimation of $g(\theta)$ based on record values does have less amount of error, since it is clearly seen in the figure that $\alpha(4) < \alpha(11)$.

## 3. Maximum Likelihood Estimations for the PDF and CDF from the Exponential Family of the First Type Based on Random Sample and Upper Record Values

Suppose that a random sample of size n from the exponential family of the first type $(X_1,...,X_n)$ has m values from the upper record $(R_1,...,R_m)$. Considering the relation (4) and the point that the maximum likelihood estimations are stable, the ML estimations for $f$ and $F$ are determined like follows:

$$\hat{\theta}_{MLE,(X_1,...,X_n)} = B^{-1}(\frac{n}{\sum_{i=1}^{n} A(X_i)}),$$

$$\hat{f}_{MLE,(X_1,...,X_n)} = A'(x)B(B^{-1}(\frac{n}{\sum_{i=1}^{n} A(X_i)}))\exp(-A(x)B(B^{-1}(\frac{n}{\sum_{i=1}^{n} A(X_i)})))$$

$$= \frac{nA'(x)}{\sum_{i=1}^{n} A(X_i)} \exp(-\frac{nA(x)}{\sum_{i=1}^{n} A(X_i)})$$

(15)

Similarly the maximum likelihood estimation for cumulative distribution function based on random sample is determined as

$$\hat{F}_{MLE,(X_1,...,X_n)} = 1 - \exp(-\frac{nA(x)}{\sum_{i=1}^{n} A(X_i)}). \quad (16)$$

In addition, the maximum likelihood estimators $f_\theta(x)$ and $F_\theta(x)$ based on the upper record values in this sample (sample of size n) are obtained as follows:

$$\hat{f}_{MLE,(R_1,...,R_m)} = \frac{mA'(x)}{A(R_m)} \exp(-\frac{mA(x)}{A(R_m)}), \hat{F}_{MLE,(R_1,...,R_m)} = 1 - \exp(-\frac{mA(x)}{A(R_m)}) \quad (17)$$

**Theorem 2.** If the random sample $(X_1,..., X_n)$ has m values from the upper records $(R_1,..., R_m)$, then:

The ML estimations of the PDF and CDF based on upper record values and based on a random sample of size m are distributed identically.

Proof: it is easily proven by considering (6), (10), (15), (16), and (17).

**Theorem 3.** The ML estimations based on random sample and based on the upper record values are biased estimators and

$$i - E[\hat{f}_{MLE,(X_1,...,X_n)}] = \sum_{i=0}^{n-2} [\frac{\Gamma(n-i-1)}{\Gamma(n)\Gamma(i+1)} (-1)^i \{nB(\theta)\}^{i+1} A'(x)]\{A(x)\}^i$$

$$ii - E[\hat{f}_{MLE,(R_1,...,R_m)}] = \sum_{i=0}^{m-2} [\frac{\Gamma(m-i-1)}{\Gamma(m)\Gamma(i+1)} (-1)^i \{mB(\theta)\}^{i+1} A'(x)]\{A(x)\}^i$$

$$iii - E[\hat{F}_{MLE,(X_1,...,X_n)}] = 1 - \sum_{i=0}^{n-1} [\frac{\Gamma(n-i)}{\Gamma(n)\Gamma(i+1)} (-1)^i \{nB(\theta)\}^i]\{A(x)\}^i$$

$$iv - E[\hat{F}_{MLE,(R_1,...,R_m)}] = 1 - \sum_{i=0}^{m-1} [\frac{\Gamma(m-i)}{\Gamma(m)\Gamma(i+1)} (-1)^i \{mB(\theta)\}^i]\{A(x)\}^i$$

Proof: Here, the forth part is proven, other parts are proven similarly:

$$\hat{F}_{MLE,(R_1,...,R_m)} = 1 - \exp(-\frac{mA(x)}{A(R_m)}), T = A(R_m) \stackrel{distribution}{=} Gamma(m, \frac{1}{B(\theta)}),$$

$$E[\hat{F}_{MLE,(R_1,...,R_m)}] = E[1 - \exp(-\frac{mA(x)}{A(R_m)})] = 1 - E[\exp(-\frac{mA(x)}{A(R_m)})]$$

$$E[\exp(-\frac{mA(x)}{A(R_m)})] = E[\sum_{i=0}^{+\infty} (\frac{(-1)^i m^i A^i(x)}{T^i i!})] = E[\sum_{even-is}^{+\Delta} (\frac{m^i A^i(x)}{T^i i!})] - E[\sum_{odd-is}^{+\infty} (\frac{m^i A^i(x)}{T^i i!})]$$

Notice that the variables $\frac{m^i A^i(x)}{T^i i!}$ are positive, that's why it is possible to interchange the expectations and summation.

$$E[\sum_{even-is}^{+\infty}(\frac{m^i A^i(x)}{T^i i!})] - E[\sum_{odd-is}^{+\infty}(\frac{m^i A^i(x)}{T^i i!})] = \sum_{even-is}^{+\infty} E(\frac{m^i A^i(x)}{T^i i!}) - \sum_{odd-is}^{+\infty} E(\frac{m^i A^i(x)}{T^i i!})$$

$$= \sum_{even-is}^{+\infty} \int_0^{+\infty} \frac{m^i A^i(x)}{t^i i!} \frac{t^{m-1} B^m(\theta) \exp(-B(\theta)t)}{\Gamma(m)} dt - \sum_{odd-is}^{+\infty} \int_0^{+\infty} \frac{m^i A^i(x)}{t^i i!} \frac{t^{m-1} B^m(\theta) \exp(-B(\theta)t)}{\Gamma(m)} dt$$

$$= \sum_{even-is}^{+\infty} \frac{m^i A^i(x) B^m(\theta)}{i! \Gamma(m)} \int_0^{+\infty} t^{m-i-1} \exp(-B(\theta)) dt - \sum_{odd-is}^{+\infty} \frac{m^i A^i(x) B^m(\theta)}{i! \Gamma(m)} \int_0^{+\infty} t^{m-i-1} \exp(-B(\theta)) dt,$$

using (14) instead of $\int_0^{+\infty} t^{m-i-1} \exp(-B(\theta)) dt$ above expression is resulted as:

$$= \sum_{even-is}^{+\infty} \frac{m^i A^i(x) B^m(\theta) \Gamma(m-i)}{i! \Gamma(m) B^{m-i}(\theta)} - \sum_{odd-is}^{+\infty} \frac{m^i A^i(x) B^m(\theta) \Gamma(m-i)}{i! \Gamma(m) B^{m-i}(\theta)} = \sum_{i=0}^{+\infty} \frac{\Gamma(m-i)}{\Gamma(i+1)\Gamma(m)} (-mA(x)B(\theta))^i,$$

finally, having considered the above relations and the fact that Gamma function is defined for positive values the expectation is readily resulted as follows:

$$E[\hat{F}_{MLE,(R_1,...,R_m)}] = 1 - E[\exp(-\frac{mA(x)}{A(R_m)})] = 1 - \sum_{i=0}^{m-1} \frac{\Gamma(m-i)}{\Gamma(i+1)\Gamma(m)} (-mA(x)B(\theta))^i.$$

**Theorem 4.**

$$i - MSE[\hat{f}_{MLE,(X_1,...,X_n)}] = \sum_{i=0}^{n-3} [(B(\theta)A^{'}(x))^2 \frac{(-2nB(\theta)A(x))^i \Gamma(n-i-2)}{\Gamma(n)\Gamma(i+1)}]$$

$$-2nA^{'}(x)(B(\theta))^2 \exp\{-A(x)B(\theta)\} \sum_{i=0}^{n-2} [\frac{(-nA(x)B(\theta))^i \Gamma(n-i-1)}{\Gamma(n)\Gamma(i+1)}] - (A^{'}(x)B(\theta))^2 \exp\{-2A(x)B(\theta)\}$$

$$ii - MSE[\hat{f}_{MLE,(R_1,...,R_m)}] = \sum_{i=0}^{m-3} [(B(\theta)A^{'}(x))^2 \frac{(-2mB(\theta)A(x))^i \Gamma(m-i-2)}{\Gamma(m)\Gamma(i+1)}]$$

$$-2mA^{'}(x)(B(\theta))^2 \exp\{-A(x)B(\theta)\} \sum_{i=0}^{m-2} [\frac{(-mA(x)B(\theta))^i \Gamma(m-i-1)}{\Gamma(m)\Gamma(i+1)}] - (A^{'}(x)B(\theta))^2 \exp\{-2A(x)B(\theta)\}$$

$$iii - MSE[\hat{F}_{MLE,(X_1,...,X_n)}] = \sum_{i=0}^{n-1} [\frac{(-2nB(\theta)A(x))^i \Gamma(n-i)}{\Gamma(n)\Gamma(i+1)}]$$

$$-2\exp\{-A(x)B(\theta)\} \sum_{i=0}^{n-1} [\frac{(-nB(\theta)A(x))^i \Gamma(n-i)}{\Gamma(n)\Gamma(i+1)}] - \exp\{-2A(x)B(\theta)\}$$

$$iv - MSE[\hat{F}_{MLE,(R_1,...,R_m)}] = \sum_{i=0}^{m-1} [\frac{(-2mB(\theta)A(x))^i \Gamma(m-i)}{\Gamma(m)\Gamma(i+1)}]$$

$$-2\exp\{-A(x)B(\theta)\} \sum_{i=0}^{m-1} [\frac{(-mB(\theta)A(x))^i \Gamma(m-i)}{\Gamma(m)\Gamma(i+1)}] - \exp\{-2A(x)B(\theta)\}$$

.

Proof: Here the forth part is proven. Other parts are proved in the same way.

$$MSE(\hat{F}_{MLE,(R_1,...,R_m)}) = E[\hat{F}_{MLE,(R_1,...,R_m)} - F_\theta(x)]^2 = E[\{1-\exp(-\frac{mA(x)}{A(R_m)})\} - \{1-\exp(-B(\theta)A(x))\}]^2 =$$

$$E[\exp(-\frac{mA(x)}{A(R_m)}) - \exp(-B(\theta)A(x))]^2 = E[\exp(-\frac{2mA(x)}{A(R_m)})] - 2\exp(-B(\theta)A(x))E[\exp(-\frac{mA(x)}{A(R_m)})] + \exp(-2B(\theta)A(x))$$

From equation number (10) it is clear that

$$T = A(R_m) \overset{distribution}{=} Gamma(n, \frac{1}{B(\theta)})$$

Considering above relation,

$$MSE(\hat{F}_{MLE,(R_1,...,R_m)}) = E(\exp(\frac{-2mA(x)}{T})) - 2\exp(-B(\theta)A(x))E(\frac{-mA(x)}{T}) + \exp(-2B(\theta)A(x)).$$

In order to obtain the above relation we should obtain the following general expectation

$$W(\alpha) = E(\exp(\frac{-\alpha mA(x)}{T})) = \int_0^{+\infty} \exp(\frac{-\alpha mA(x)}{T}) \frac{t^{m-1}B^m(\theta)\exp(-B(\theta)t)}{\Gamma(m)} =$$

$$\int_0^{+\infty} \sum_{i=0}^{+\infty} \frac{(-\alpha mA(x))^i}{t^i i!} \frac{t^{m-1}B^m(\theta)\exp(-B(\theta)t)}{\Gamma(m)} = \int_0^{+\infty} \sum_{even-indices}^{+\infty} \frac{(\alpha mA(x))^i}{t^i i!} \frac{t^{m-1}B^m(\theta)\exp(-B(\theta)t)}{\Gamma(m)} -$$

$$\int_0^{+\infty} \sum_{odd-indices}^{+\infty} \frac{(\alpha mA(x))^i}{t^i i!} \frac{t^{m-1}B^m(\theta)\exp(-B(\theta)t)}{\Gamma(m)}.$$

It is clear that $\frac{(\alpha mA(x))^i}{t^i i!} > 0$ then,

$$\int_0^{+\infty} \sum_{even-indices} \frac{(\alpha mA(x))^i}{t^i i!} \frac{t^{m-1}B^n(\theta)\exp(-B(\theta)t)}{\Gamma(m)} dt - \int_0^{+\infty} \sum_{odd-indices} \frac{(\alpha mA(x))^i}{t^i i!} \frac{t^{m-1}B^m(\theta)\exp(-B(\theta)t)}{\Gamma(m)} dt =$$

$$\sum_{even-indices} \int_0^{+\infty} \frac{(\alpha mA(x))^i}{t^i i!} \frac{t^{m-1}B^m(\theta)\exp(-B(\theta)t)}{\Gamma(m)} dt - \sum_{odd-indices} \int_0^{+\infty} \frac{(\alpha mA(x))^i}{t^i i!} \frac{t^{m-1}B^m(\theta)\exp(-B(\theta)t)}{\Gamma(m)} dt =$$

$$\sum_{even-indices} \frac{(\alpha mA(x))^i B^m(\theta)}{\Gamma(i+1)\Gamma(m)} \int_0^{+\infty} \frac{t^{m-1}\exp(-B(\theta)t)}{t^i} dt - \sum_{odd-indices} \frac{(\alpha mA(x))^i B^m(\theta)}{\Gamma(i+1)\Gamma(m)} \int_0^{+\infty} \frac{t^{m-1}\exp(-B(\theta)t)}{t^i} dt =$$

$$\sum_{even-indices} \frac{(\alpha mA(x))^i B^m(\theta)}{\Gamma(i+1)\Gamma(m)} \int_0^{+\infty} t^{m-i-1}\exp(-B(\theta)t)dt - \sum_{odd-indices} \frac{(\alpha mA(x))^i B^m(\theta)}{\Gamma(i+1)\Gamma(m)} \int_0^{+\infty} t^{m-i-1}\exp(-B(\theta)t)dt$$

$$\sum_{even-indices} \frac{(\alpha mA(x))^i B^m(\theta)}{\Gamma(i+1)\Gamma(m)} \frac{\Gamma(m-i)}{B^{m-i}(\theta)} - \sum_{odd-indices} \frac{(\alpha mA(x))^i B^m(\theta)}{\Gamma(i+1)\Gamma(m)} \frac{\Gamma(m-i)}{B^{m-i}(\theta)} =$$

$$\sum_{i=0}^{+\infty} \frac{(\alpha mA(x))^i B^m(\theta)}{\Gamma(i+1)\Gamma(m)} \frac{\Gamma(m-i)}{B^{m-i}(\theta)} = \sum_{i=0}^{+\infty} \frac{(-\alpha mB(\theta)A(x))^i \Gamma(m-i)}{\Gamma(i+1)\Gamma(m)} = \sum_{i=0}^{n-1} \frac{(-\alpha mB(\theta)A(x))^i \Gamma(m-i)}{\Gamma(i+1)\Gamma(m)}.$$

Finally, $W(\alpha) = \sum_{i=0}^{m-1} \frac{(-\alpha mB(\theta)A(x))^i \Gamma(m-i)}{\Gamma(i+1)\Gamma(m)}$.

Now the MSE function is obtained easily as follow:

$$MSE = W(2) - 2\exp(-B(\theta)A(x))W(1) + \exp(-2B(\theta)A(x)) =$$

$$\sum_{i=0}^{m-1} \frac{(-2mB(\theta)A(x))^i \Gamma(m-i)}{\Gamma(i+1)\Gamma(m)} - 2\exp(-B(\theta)A(x))\sum_{i=0}^{m-1} \frac{(-mB(\theta)A(x))^i \Gamma(m-i)}{\Gamma(i+1)\Gamma(m)} + \exp(-2B(\theta)A(x)).$$

The next theorem clarifies some advantageous aspects of the above likelihood estimators. That is, the above-mentioned estimators are always nearly unbiased in large samples.

**Theorem 5.** The maximum likelihood estimations based on random sample and on the upper record values for the PDF and CDF in the exponential family of the first type are asymptotically unbiased. In the other words

$i - \underset{n\to\infty}{Lim} E[\hat{f}_{MLE,(X_1,...,X_n)}] = f_\theta(x)$

$ii - \underset{m\to\infty}{Lim} E[\hat{f}_{MLE,(R_1,...,R_m)}] = f_\theta(x)$

$iii - \underset{n\to\infty}{Lim} E[\hat{F}_{MLE,(X_1,...,X_n)}] = F_\theta(x)$

$iv - \underset{m\to\infty}{Lim} E[\hat{F}_{MLE,(R_1,...,R_m)}] = F_\theta(x)$

The following lemma is needed in order to prove the above theorem.

Lemma 1. $\underset{n\to\infty}{Lim} \frac{\Gamma(n-i-1)n^{i+1}}{\Gamma(n)} = 1$ for the reason that

$$\underset{n\to\infty}{Lim}\frac{\Gamma(n-i-1)n^{i+1}}{\Gamma(n)} = \underset{n\to\infty}{Lim}\frac{n^i}{\prod_{k=1}^{i+1}(n-k)} = \underset{n\to\infty}{Lim}\frac{n^i}{n^i} = 1$$

It must be proven that

$$\forall \in > 0, \exists N, s.t, \forall n \geq N : \frac{\Gamma(n-i-1)n^{i+1}}{\Gamma(n)} - 1 < \in$$

Now, assuming that $\in$ is arbitrary, N is chosen as $N = \frac{\sqrt[i+1]{(\in+1)} - i\sqrt[i+1]{(\in+1)}}{\sqrt[i+1]{(\in+1)} - 1}$. Now

$$\frac{\Gamma(n-i-1)n^{i+1}}{\Gamma(n)} = \frac{(n-i-2)!n^{i+1}}{(n-1)!} = \frac{n^{i+1}}{(n-1)(n-2)...(n-i-1)} = \underbrace{\frac{n}{n-1}\times...\times\frac{n}{n-i-1}}_{(i+1)times} \leq \underbrace{\frac{n}{n-i-1}\times...\times\frac{n}{n-i-1}}_{(i+1)times}$$

On the other hand, it can be simply shown that

$$\forall n \geq N, and, i \in \{1,2,3,...,N-2\} : \frac{n}{n-i-1} \leq \frac{N}{N-i-1},$$

thus

$$\underbrace{\frac{n}{n-i-1}\times...\times\frac{n}{n-i-1}}_{(i+1)times} \leq \underbrace{\frac{N}{N-i-1}\times...\times\frac{N}{N-i-1}}_{(i+1)times},$$

and consequently

$$\forall n \geq N : \frac{\Gamma(n-i-1)n^{i+1}}{\Gamma(n)} < \underbrace{\frac{N}{N-i-1}\times...\times\frac{N}{N-i-1}}_{(i+1)times}.$$

Finally by substituting $N = \frac{\sqrt[i+1]{(\in+1)} - i\sqrt[i+1]{(\in+1)}}{\sqrt[i+1]{(\in+1)} - 1}$ in the last relation, the following result is reached:

$$\frac{\Gamma(n-i-1)n^{i+1}}{\Gamma(n)} < \in +1$$

**Proof:**

Only the first part is proven; the other parts have similar proofs. First notice that

$$Lim_{n\to\infty}[\hat{f}_{MLE,(X_1,...,X_n)}] = \sum_{i=0}^{\infty}[Lim_{n\to\infty}\frac{\Gamma(n-i-1)}{\Gamma(n)\Gamma(i+1)}(-1)^i n^{i+1}\{B(\theta)\}^{i+1}A'(x)]\{A(x)\}^i]$$

$$= \sum_{i=0}^{\infty} Lim_{n\to\infty}(\frac{\Gamma(n-i-1)n^{i+1}}{\Gamma(n)}\frac{\{-B(\theta)\}^i B(\theta)A'(x)\{A(x)\}^i}{\Gamma(i+1)}]$$

By considering lemma (1)

$$Lim_{n\to\infty}[\hat{f}_{MLE,(X_1,...,X_n)}] = \sum_{i=0}^{\infty}\frac{\{-A(x)B(\theta)\}^i}{i!}B(\theta)A'(x) = A'(x)B(\theta)\exp\{-A(x)B(\theta)\}.$$

The above theorem makes us to pay attention to convergence in probability and to the following theorem.

**Theorem 6.** If a random sample $(X_1,..., X_n)$ has m values from the upper records $(R_1,..., R_m)$, then the maximum likelihood estimations for the CDF and PDF based on upper record values are consistent estimations for the CDF and PDF.

$$i - \hat{f}_{MLE,(R_1,...,R_m)} \xrightarrow{P} f_\theta(x)$$

$$ii - \hat{F}_{MLE,(R_1,...,R_m)} \xrightarrow{P} F_\theta(x)$$

**Proof.** The part i is proven; other part has a similar proof. It is considered the following abbreviation just for easiness.

$$\hat{f}_{MLE,(R_1,...,R_n)} = \hat{f}_{MLE,R}$$

Based on the Markov's theorem

$$P\{\varepsilon \le |\hat{f}_{MLE,R} - E[\hat{f}_{MLE,R}]|\} \le \frac{Var(\hat{f}_{MLE,R})}{\varepsilon}$$

On the other hand, it is concluded by some mathematical calculations that

$$E[\{\hat{f}_{MLE,R}\}^2] = \sum_{i=0}^{+\infty}\frac{\Gamma(m-i-2)}{\Gamma(i+1)\Gamma(m)}\{A'(x)\}^2\{B(\theta)\}^2(-2)^i\{nB(\theta)\}^i\{A(x)\}^i$$

and consequently

$$Var(\hat{f}_{MLE,R}) = \sum_{i=0}^{+\infty}\frac{\Gamma(m-i-2)}{\Gamma(i+1)\Gamma(m)}\{A'(x)\}^2\{B(\theta)\}^2(-2)^i\{mB(\theta)\}^i\{A(x)\}^i - \{E[\hat{f}_{MLE,R}]\}^2$$

thus

$$Lim_{m\to\infty} Var(\hat{f}_{MLE,R}) = \sum_{i=0}^{+\infty} Lim_{m\to\infty}[\frac{m^i \Gamma(m-i-2)}{\Gamma(m)}] \frac{\{-2B(\theta)\}^i \{A(x)\}^i}{\Gamma(i+1)} \{A'(x)\}^2 \{B(\theta)\}^2 - \{f_\theta(x)\}^2$$

$$= \{A'(x)\}^2 \{B(\theta)\}^2 \exp\{-2A(x)B(\theta)\} - \{f_\theta(x)\}^2 = 0$$

Therefore $\quad Lim_{m\to\infty} P\{|\hat{f}_{MLE,R} - E[\hat{f}_{MLE,R}]| > \varepsilon\} \leq 0 \quad$ so $\quad \hat{f}_{MLE,R} \xrightarrow{P} E[\hat{f}_{MLE,R}] \quad$ or

$$\hat{f}_{MLE,R} - E[\hat{f}_{MLE,R}] \xrightarrow{P} 0$$

On the other hand, it is reminded that if $Lim_{m\to\infty} a_m = a$ and $X \xrightarrow{P} r$ then (Billingsley (1995))

$$X + a_m \xrightarrow{P} r + a \tag{18}$$

Considering (18) and respecting the second part of theorem 5 ($Lim_{m\to\infty} E[\hat{f}_{MLE,R}] = f_\theta(x)$), it is obtained that

$$\hat{f}_{MLE,R} - E[\hat{f}_{MLE,R}] + E[\hat{f}_{MLE,R}] - f \xrightarrow{P} 0,$$

or equally

$$\hat{f}_{MLE,R} \xrightarrow{P} f.$$

The following remembrance is essential before expressing the next theorem.

Remembrance 1: if $\{X_n\} \xrightarrow{P} \{X\}$ and $D$ is a continuous function of $y$ $\{D(y): R \to R\}$, then (Billingsley, 1995)

$$D(\{X_n\}) \xrightarrow{P} D(\{X\})$$

**Theorem 7.**

Suppose that a random sample of size n from the exponential family of the first type is available. Also suppose that this sample includes m values from the upper record values. The estimators $\hat{\theta}_{MLE,(R_1,\ldots,R_m)}$ are consistent for $\theta$. Namely

$$\hat{\theta}_{MLE,(R_1,\ldots,R_m)} \xrightarrow{P} \theta$$

**Proof.**

With respect to the continuity of functions $B(t), -Log(1-t), F(t) = 1 - \exp\{-A(t)B(\theta)\}$, and the continuity theorem of composite functions, it is concluded that the function $D(t) = B(\frac{-Log(1-F(t))}{A(x)})$ is always continuous on $\mathbb{R}$. On the other hand, it is known from the second part of theorem 6 that

$$\hat{F}_{MLE,(R_1,...,R_m)} \xrightarrow{P} F_\theta(x),$$

or equivalently

$$A(x)B(\hat{\theta}_{MLE,(R_1,...,R_m)})\exp\{-A(x)B(\hat{\theta}_{MLE,(R_1,...,R_m)})\} \xrightarrow{P} F_\theta(x).$$

Now with respect to remembrance 1,

$$D(\hat{F}_{MLE,(R_1,...,R_m)}) \xrightarrow{P} D(F_\theta(x))$$

which is equivalent to

$$\hat{\theta}_{MLE,(R_1,...,R_m)} \xrightarrow{P} \theta$$

## 4. Conclusions

In the beginning, considering two general subfamilies from the exponential family, it was tried to determine the maximum likelihood estimation for these two families at the first step and general state; then, the maximum likelihood estimations based on upper were calculated for the unknown parameter of these two families; it was engaged via some theorems to description of this problem that the record estimations can be acceptable in general state; then the estimations of maximum likelihood and of maximum likelihood based on records were calculated for the cumulative distribution function and for the probability density function of these families, and their corresponding errors were described and studied by using some theorems; in addition, we engaged to proof of this problem that at asymptotical states, these estimations are unbiased estimators. Subsequently, several theorems associated with convergence in probability were expressed for the resulted estimators.